\documentclass[11pt,leqno]{article}
\usepackage{epsfig}
\usepackage{graphicx}
\usepackage{wrapfig}
\usepackage{amsmath}
\usepackage{amssymb}
\usepackage{amscd }
\usepackage{amsthm}

\usepackage{tikz}
\usetikzlibrary{arrows}

\newcommand{\be}{\begin{eqnarray}}
\newcommand{\ee}{\end{eqnarray}}
\newcommand{\bi}{\begin{itemize}}
\newcommand{\ei}{\end{itemize}}
\newcommand{\bd}{\begin{definition}}
\newcommand{\ed}{\end{definition}}
\newcommand{\bt}{\begin{theorem}}
\newcommand{\et}{\end{theorem}}
\newcommand{\bc}{\begin{corollary}}
\newcommand{\ec}{\end{corollary}}
\newcommand{\bcn}{\begin{conjecture}}
\newcommand{\ecn}{\end{conjecture}}
\newcommand{\br}{\begin{remark}}
\newcommand{\er}{\end{remark}}
\newcommand{\ce}{\begin{eqnarray*}}
\newcommand{\de}{\end{eqnarray*}}
\newcommand{\bpf}{\begin{proof}}
\newcommand{\epf}{\end{proof}}
\newcommand{\bl}{\begin{lemma}}
\newcommand{\el}{\end{lemma}}
\newtheorem{theorem}{Theorem}[section]
\newtheorem{case}{Case}
\newtheorem{conjecture}[theorem]{Conjecture}
\newtheorem{lemma}[theorem]{Lemma}
\newtheorem{remark}[theorem]{Remark}
\newtheorem{definition}[theorem]{Definition}
\newtheorem{proposition}[theorem]{Proposition}
\newtheorem{Examples}[theorem]{Examples}
\newtheorem{corollary}[theorem]{Corollary}
\numberwithin{equation}{section}

\def\e{\varepsilon}

\def\a{\alpha}
\def\o{\omega}
\def\b{\beta}
\def\d{\delta}

\def\g{\gamma}
\def\l{\lambda}

\def\[{{\Big[}}
\def\]{{\Big]}}
\def\<{{\langle}}
\def\>{{\rangle}}

\def\({{\Big(}}
\def\){{\Big)}}

\def\bt{\begin{theorem}}
\def\et{\end{theorem}}
\def\bl{\begin{lemma}}
\def\el{\end{lemma}}
\def\br{\begin{remark}}
\def\er{\end{remark}}
\def\bx{\begin{Examples}}
\def\ex{\end{Examples}}
\def\bd{\begin{definition}}
\def\ed{\end{definition}}
\def\bp{\begin{proposition}}
\def\ep{\end{proposition}}
\def\bc{\begin{corollary}}
\def\ec{\end{corollary}}

\def\cB{{\mathcal B}}
\def\cC{{\mathcal C}}

\def\cG{{\mathcal G}}

\def\cR{{\mathcal R}}

\def\cT{{\mathcal T}}

\def\cZ{{\mathcal Z}}

\def\mN{{\mathbb N}}

\def\mP{{\mathbb P}}

\def\mX{{\mathbb X}}

\def\mZ{{\mathbb Z}}

\usepackage{graphicx}

\begin{document}

\section* {\center{Whirly 3-Interval Exchange Transformations}}
\begin{center}
\textsc
{
Yue Wu\footnote{Schlumberger PTS Full Waveform Inversion Center of Excellence, Houston, Texas, USA}\\
3519 Heartland Key LN, Katy, TX, 77494\\
wuuyue@gmail.com
}
\end{center}

\begin{abstract}
Irreducible interval exchange transformations are studied with regard to whirly property, a condition for non-trivial spatial factor. Uniformly whirly transformation is defined and to be further studied. An equivalent condition is introduced for whirly transformation. We will prove that almost all 3-interval exchange transformations are whirly, using a combinatorics approach with application of the Rauzy-Veech Induction. It is still an open question whether whirly property is a generic property for $m$-interval exchange transformations ($m\geq 4$).
\end{abstract}
\section* {}
\indent \indent Interval Exchange Transformations, as a set of important dynamical systems, have been actively studied for decades. We recall some of the key theorems, either the results or the methods of which are related to the current study or possible extensions of this paper in the future. The proof of the unique ergodic property of measure theoretical generic interval exchange transformations was achieved by H. Masur\cite{MASUR} and W.A. Veech\cite{VEE1} independently using geometric methods, and was proved later using mainly combinatorial methods by M. Boshernitzan\cite{BOS}. A. Avila and G. Forni\cite{AVILA-FORNI} showed that weak mixing is  a measure theoretical generic property for irreducible $m$-interval exchange transformations ($m\geq 3$). J. Chaika\cite{CHAI} developed a general result showing that any ergodic transformation is disjoint with almost all interval exchange transformations.  J. Chaika and J. Fickenscher\cite{CHAI2} showed that topological mixing is a topologically residual property for interval exchange transformations.\\
\indent The concept of whirly transformation was introduced and studied in E. Glasner, B. Weiss\cite{GW}, E. Glasner, B. Tsirelson, B. Weiss\cite{GTW}. In E. Glasner, B. Weiss\cite{GW2}, Proposition 1.9. states that the near action (weak closure of all the powers) of a transformation admits no non-trivial spatial factors if and only if it is whirly. In the Section of Introduction, we recall the relevant notions and facts about interval exchange transformations and whirly transformation. A new notion introduced in this section is that of uniformly whirly transformation. In the second section we study the space of three interval exchange transformations and deduce facts about the visitation times of Rauzy-Veech induction.
In the last section we complete the proof of the major theorem, which states almost all three interval exchange transformations are whirly, thus admit no nontrivial spatial factor. To prove this main theorem, first we establish an equivalent definition of whirly transformations based on the assumption of ergodicity. Then the major key facts are deduced as Claim 1, Claim 2 and Claim 3, which show the whirly property for the base of the Rohlin tower associated with the Veech-induction map. Finally we apply a density point argument to extend the property to arbitrary general non-null measurable sets.
\section{Introduction}\label{section1}
\indent \indent An interval exchange transformation perturbs the half-closed half-open subintervals of a  half-closed half-open interval. The subintervals have lengths corresponding to the vector
$\lambda=(\lambda_{1},\cdots ,\lambda _{m})$, $ \lambda _{i}
>0,\, 1\leq i \leq m$. All such vectors form a positive cone $\Lambda _{m}\subset
R^{m}$. The subintervals thus are
$[\beta_{i-1},\beta_{i})$,$\,1\leq i \leq m$, with
$\bigcup[\beta_{i-1},\beta_{i})=[0,\left|\l\right|)$, where\\
 \begin{equation}\begin{array}{l}
 \mid \l \mid=\sum\limits_{i=1}^{m} \l_{i}\\
 \text{and, }\\
     \beta_{i}(\lambda)=\left\{\begin{array}{clcr} 0& i=0\\
         \sum\limits_{j=1}^{i}\l_{j}& 1\leq i\leq m .
                      \end{array}\right.
 \end{array}\end{equation}
\indent Let $\cG_{m}$ be the group of m-permutations, and
$\mathcal{G}_{m}^{0}$ be the subset of $\mathcal{G}_{m}$ which
contains all the irreducible permutations on $\{1,2,\cdots, m\}$. A
permutation $\pi$ is irreducible if and only if for any $1\leq k
<m,\,\{1,2,\cdots,k\}\neq\{\pi (1),\cdots,\pi(k)\}$, or
equivalently $\sum\limits_{j=1}^{k}(\pi(j)-j)>0,\, (1\leq k<m))$. Given $\lambda \in \Lambda_{m},\, \pi \in \mathcal{G}_{m}^{0}$,
the corresponding interval exchange transformation is defined by:\\
 \begin{equation}\begin{array}{l}
 T_{\l,\pi}(x)=x-\beta_{i-1}(x)+\beta_{\pi i-1}(\lambda^{\pi}),\quad (x\in
[\beta_{i-1}(\l),\quad
\beta_{i}(\l))\,),\\
\text{where }\lambda^{\pi}=(\l_{\pi^{-1}1},\l_{\pi^{-1}2},\cdots,
\l_{\pi^{-1}m}) .
\end{array}\end{equation}
Obviously $\beta_{\pi i-1}(\l^{\pi})=\sum\limits_{j=1}^{\pi
i-1}\l_{\pi^{-1}j}$,
 and the transformation $T_{\l,\pi}$, which is also denoted by $(\l,\pi)$, sends the
$i$th interval to the $\pi (i)$th position.\\

In M. Keane\cite{KEA}, the i.d.o.c.(infinite distinct orbits condition) is raised for
the sufficient condition of minimality: $\lambda ,\pi$ is said to
satisfy the {\em i.d.o.c\/}.
if
\begin{enumerate}
\item[\em i)\/] for any $0\leq i <m, \{T^{k}\b_{i}, k\in \mZ\}$ is a infinite set;
\item[\em ii)\/] $\{T^{k}\b_{i},k\in\mZ\}\cap\{T^{k}\b_{j},k\in \mZ\}=\emptyset$, whenever $i\neq j.$
\end{enumerate}
\indent Suppose $m>1$, $(\l ,\pi)\in\Lambda _{m}\times \cG ^{*}_{m}$ , where $\cG ^{*}_{m}$ is the set
of irreducible permutations with the property that $\pi (j+1)\neq
\pi (j)+1$ for all $1\leq j \leq m-1$. Let $I$ be an interval of
the form $I=[\xi ,\eta)$, $0\leq \xi <\eta \leq \left|\l\right|$.
Since $T$ is defined on $[0,\left|\lambda\right|)$, and $T$ is Lebesgue measure
preserving, we know that Lebesgue almost all points of $I$ return
to $I$ infinitely often under iteration of $T$. We use $T|_{I}$ to
denote the induced transformation of $T$
 on $I$. By W.A. Veech\cite{VEE4}, $T|_{I}$ is an interval exchange transformation with $(m-2)$, $(m-1)$, or $m$ discontinuities.
\bd [Admissible Interval; W.A. Veech \cite{VEE4}] \label{AdmInt}
Suppose $(\l ,\pi)$ satisfies the i.d.o.c., and $I=(\xi
,\eta)$ where $\xi =T^{k}\b _{s}$,$(1\leq s<m)$; $\eta =T^{l}\b
_{t}$ $(1\leq t<m)$, and $\tau \in\{k,l\}$ have the following
property: If $\tau \geq 0$, there is no $j$, $0<j<\tau$, such that
$T^{j}\b_{s}\in I$; If $\tau< 0$, there is no $j$, $0\geq j>\tau$,
such that
$T^{j}\b_{s}\in I$. Then we say that $I$ is an admissible subinterval of $(\l ,\pi)$.
\ed
{\em\bf Rauzy-Veech induction.\/\/} For $T_{\l,\pi}$, the Rauzy map
sends it to its induced map on
$[0,\left|\l\right|-min\left\{\l_{m},\l_{\pi^{-1}m}\right\})$,
which is the largest admissible interval of form
$J=[0,L),0<L<\left|\l\right|$.\\
Given any permutation, two actions $a$ and $b$ are:
\be
a(\pi)(i)=\left\{\begin{array}{llcr}\pi(i) & i\leq\pi^{-1}m
\\\pi(i-1)&\pi^{-1}m+1<i\leq m \\ \pi(m)
&i=\pi^{-1}m+1 \end{array}\right.
\ee
 and
 \be
b(\pi)(i)=\left\{\begin{array}{llcr}\pi(i) & \pi(i)\leq\pi(m)
\\\pi(i)+1&\pi(m)+1<\pi(i)< m \\ \pi(m)+1
&\pi(i)=m . \end{array}\right.
\ee
 \indent The Rauzy-Veech map $\mathcal{Z}(\l,\pi):\,
\Lambda_{m}\times \mathcal{G}_{m}^{0}\rightarrow \Lambda_{m}
\times \mathcal{G}_{m}^{0}$ is determined by :
\be\label{RzVch}\mathcal{Z}(\l,\pi)=(A(\pi ,c)^{-1}\l,c\pi) ,
\ee
where $c=c(\l,\pi)$ is defined by
\be
c(\l,\pi)=\left\{\begin{array}{clcr}a,&\l_{m}<\l_{\pi^{-1}m}\\
b,&\l_{m}>\l_{\pi^{-1}m.}\end{array}\right.
\ee
$\mathcal{Z}(\l,\pi)$ is a.e. defined on
$\Lambda_{m}\times\left\{\pi\right\}$, for each $\pi\in\mathcal{G}_{m}^{0}$. \\
\indent The matrices $A=A(\pi,c)$ in \ref{RzVch} are defined as the following:
\begin{equation} A(\pi ,a)=\left(\begin{tabular}{c|c}
$I_{\pi^{-1}m}$&$\begin{array}{ccccc}
0&0&\cdots&0&0\\
0&0&\cdots&0&0\\
.&.&\cdots&.&.\\
0&0&\cdots&0&0\\
1&0&\cdots&0&0
  \end{array}$\\\hline\\
0& $\begin{array}{ccccc}
0&1&\cdots&0&1\\
0&0&\cdots&0&0\\
.&.&\cdots&.&.\\
0&0&\cdots&0&1\\
1&0&\cdots&1&0
  \end{array}$\\
\end{tabular}\right)
\end{equation}
\be A(\pi ,b)=\left(
\begin{tabular}{c|c}
 $I_{m-1}$ &0\\\hline\\
$\underbrace{\begin{array}{ccccccc}
0&\cdots&0&1&0&\cdots&0\end{array}}_{\mbox{1 at the jth position}}$&1\\
\end{tabular}
\right)
\ee
where $I_{k}$ is the $k$-identity matrix, and $j=\pi^{-1}m$.\\
\indent And the normalized Rauzy map $\cR : \, \Delta_{m-1}\times \cG^{0}_{m}\rightarrow
 \Delta_{m-1}\times\cG^{0}_{m}$  is defined by
 \be
 \cR (\l, \pi)=(\frac{\displaystyle A(\pi ,c)^{-1}\l}{\displaystyle \left|A(\pi ,c)^{-1}\l\right|},
 c\pi)=(\frac{\displaystyle \pi^{*}_{1}\cZ (\l
,\pi)}{\displaystyle \left|\pi^{*}_{1}\cZ (\l ,\pi)\right|},
\pi^{*}_{2}\cZ (\l ,\pi)),
 \ee
 where $\pi^{*}_{1}$ and $\pi^{*}_{2}$ are the projection to the first coordinate and the second coordinate respectively.\\
Iteratively,
\be\label{Cn}
\mathcal{Z}^{n}(\l,\pi)=((A^{(n)})^{-1}\l,
c^{(n)}\pi)=(\l ^{(n)},\pi^{(n)}),
\ee
where
\be\label{An}
c^{(n)}=c_{n}c_{n-1}\cdots c_{1},(c_{1},\cdots,c_{n}\in
 \{a,b\}, c_{i}=c(\mathcal{Z} ^{i-1}(\l,\pi)))
\ee
and
\be
A^{(n)}=A(\pi,c_{1})A(c^{(1)}\pi, c_{2})A(c^{(2)}\pi,c_{3})\cdots
 A(c^{(n-1)}\pi,c_{n}) .
 \ee
\indent The Rauzy class $\cC\subseteq\cG_{m}$ of $\pi$ is a set of orbits for the group of maps generated by $a$ and $b$. On the $\cR$ invariant component $\Delta_{m-1}\times\cC$, we have:
\bt [H.Masur\cite{MASUR};W.A. Veech\cite{VEE1}]\label{Ergodic} Let $\pi\in\cG_{m}^{0}$, the set of irreducible permutations. For Lebesgue almost all $\l\in\Lambda_{m}$, normalized Lebesgue measure on $I^{\l}$ is the unique invariant Borel probability measure for $T_{(\l,\pi)}$. In particular, $T_{(\l,\pi)}$ is ergodic for almost all $\l$.
\et\
{\em\bf Whirly Action, Whirly Automorphism. \/\/}
In this paper, we assume weak topology as defined in the following Definition \ref{WeakTp}.
\bd [Weak Topology on Automorphism Group] \label{WeakTp} Let $(\mX ,\cB ,\mu)$ be a standard probability Borel space, and $G=Aut(\mX ,\cB ,\mu)$ be the group of all non-singular measurable automorphisms of $(\mX ,\cB ,\mu)$. Suppose $(E_{n})$ is a countable family of measurable subsets generating  $\cB$. The weak topology of $G$ is generated by the metric $d(S,T)$, for any $S,\,T \in G$, where $d(S,T)=\sum_{n=1}^{\infty}2^{-n}\mu(SE_{n}\triangle TE_{n})$. \ed
Utilizing the weak topology defined as above, one can project the concept of whirly action (Definition \ref{WhirlyAc}) to whirly automorphism (Definition \ref{whirly1}). This is included in the following review of the definitions and fundamental propositions about whirly action and whirly automorphisms.\\

Whirly action is introduced by E. Glasner, B. Tsirelson, B. Weiss \cite{GTW}, Definition 3.1. The purpose is to study the condition for a Polish group action to admit a spatial model. In the same paper, they translated the concept of whirly from the group action to automorphisms since the weak closure of a rigid automorphism is a near action. They showed that in the
group $G$ of automorphisms on a finite Lebesgue space, whirly (in
the sense of $Z$ action) is a topologically generic property, i.e. the
set of whirly automorphisms is residual in $G$. The concept of `whirly transformation' is inherited from the theory
about general group actions, and implies weak mixing.  It is interesting to
ask whether whirly is a generic property in the space of interval
exchange transformations. Theorem \ref{Main} gives a positive answer for three interval exchange transformations.\\
\indent Without considering the measure, we have the Borel action,
satisfying similar condition as in Definition ¡¢\ref{Near}, defined
below:
\bd [Borel Action]\label{Borel}
Suppose $G$ is a Polish group and $(\mX ,\cB ,\mu)$ is a standard probability Borel
space. We say a Borel map $G \times \mX \rightarrow \mX $ $((g,
x)\rightarrow gx)$ is a Borel action of $G$ on $(\mX ,\cB ,\mu)$ if
it satisfies the following properties:\\
{\em (i)\/} $\quad$ $ex=x$ for all $x\in \mX$, where $e$ is the
identity element of $G$;\\
{\em (ii)\/} $\quad$ $g(hx)=(ghx)$ for all $x\in \mX$, where $g,
 h\in G$.
 \ed
  \bd [Spatial $G$ Action: E. Glasner, B. Tsirelson, B. Weiss\cite{GTW}]\label{Spatial}
  A spatial $G$-action on a standard Lebesgue space $(\mX ,\cB
  ,\mu)$ is
 a Borel action of $\mP$ on the space such that each $g\in \mP$
 preserves $\mu$.
\ed
The concept of near action is introduced measure theoretically:
\bd[Near Action: E. Glasner, B. Tsirelson, B. Weiss\cite{GTW}]\label{Near}
Suppose $\mP$ is a Polish group and $(\mX ,\cB ,\mu)$ is a
standard probability Borel space. We say  a Borel map $\mP \times
\mX \rightarrow \mX $ $((g, x)\rightarrow gx)$ is a near action of
$\mP$ on $(\mX ,\cB ,\mu)$ if it satisfies the following
properties:\\
{\em (i)\/} $\quad$ $ex=x$ for a.e. $x\in \mX$, where $e$ is the
identity element of $\mP$;\\
{\em (ii)\/} $\quad$ $g(hx)=(ghx)$ for a.e. $x\in \mX$, where $g, h\in\mP$;\\
{\em (iii)\/} $\quad$Each $g\in G$ preserves the measure $\mu$.
\ed
   {\bf Note.\/} the set of measure one in Definition \ref{Near} (ii) may depend on the
   pair $g,h$. It is easy to see that a near action is a continuous
   homomorphism from $\mP$ to $G$ ($G$ is the automorphism group of $\mX$).\\

Now we define the key concept of this paper:
\bd
[Whirly Action: E. Glasner, B. Tsirelson, B. Weiss\cite{GTW}]\label{WhirlyAc}
$\,$Given $\e>0$, if for all sets $E,F\in \cB$ with $\mu (E),
\mu (F)>0$, there exists $\g\in N_{\e} (Id)$ (the $\e$ neighborhood
of the identity $Id=e$ in $\mP$), such that $\mu (E\cap \g F)>0$ then we
say the near action of $\mP$ on $(\mX ,\cB
,\mu)$ is whirly.
\ed
\bt[E. Glasner, B. Tsirelson, B. Weiss\cite{GTW} Proposition 3.3]\label{Factor}
A whirly action
does not admit a nontrivial spatial factor, and thus has no spatial model.
\et
\br
If an automorphism $(\mX ,\cB ,T
,\mu)$ is rigid, then its weak closure '$Wcl(T)$' is a closed
subgroup of $G=Aut(\mX ,\cB ,\mu)$. With the induced topology,
$Wcl(T)$ is also a Polish space. Based on this fact, the whirly
transformation is a concept induced from whirly action.
\er
Let $(\mX,\cB, \mu )$ be the standard Lebesgue probability
space, $\mX=[0,1]$, and denote $G=Aut(\mX)$ the Polish group of its automporphism.
\bd[Whirly Automorphism]\label{whirly1}
We say a rigid system $(X,\cB,\mu,
T)$ is whirly, if given $\epsilon >0$ for any $\mu$ positive measure
sets $E$ and $F$ $(\mu (E), \mu (F)>0)$ in $\cB$, there exists $n$
such that $T^{n}\in U_{\epsilon}$  (the
$\epsilon$-neighborhood
of the identity map in the weak topology of $G$), and $\mu(T^{n}E\cap F)>0$.
\ed
Whirly implies rigid. E.Glasner, B.Weiss\cite{GW} Corollary 4.2. showed that if $(\mX,\cB,\mu,T)$ is whirly then it is weak mixing. In the same paper as Theorem 5.2., it is proved that:
\bt[E.Glasner, B.Weiss\cite{GW} Theorem 5.2]
The set of all the whirly transformations is residual (dense $G_{\d}$ subset) in $G$.
\et
Next we introduce a new notion of uniformly whirly, which is stronger than or equivalent to whirly:
\bd[Uniformly Whirly]
A rigid system $(\mX ,\cB ,\mu
,T)$ is uniformly whirly if given $\e>0$ for any $0<\a ,\b<1$, we have
$$\underset{\mu (E)=\a ,\mu (F)=\b}{inf}\quad \underset{T^{n}\in U_{\e}}{sup}\{\mu(T^{n}E \cap F)\}>0 .$$
\ed
Uniformly whirly implies whirly.\\

 \textbf{Questions:} Is uniformly whirly equivalent to whirly? If not, is the collection of uniformly whirly automorphisims a dense $G_{\delta}$ subset of $G$.\\

\indent It is interesting to
ask whether whirly property is a generic property in the space of $m$ interval
exchange transformations ($m\geq 3$). The major theorem (Theorem \ref{Main}) of this paper provides a positive answer for $m=3$. Below is the main result and an outline of the proof:
\bt\label{Main}
Let $\pi =(3,2,1)$, for Lebesgue
almost all $\l\in\Lambda_{3}$. The three dimensional cone of positive real numbers, the interval exchange transformation $T_{(\l ,\pi)}$ is whirly.
\et

{\em\bf Outline of the proof of Theorem
\ref{Main}\/\/}\\
\indent First, we raise an equivalent definition for whirly transformation (Definition \ref{whirly2}). This definition enables us to use the cyclic approximation of rank 1 stacking structure (\cite{VEE3} Section 3) associted with Rauzy-Veech induction more effectively.\\
\indent Second, for symmetric $3$-permutation $\pi$, we study the Veech induction map $\cT_{2}:\l\rightarrow
\frac{\a}{\left|\a\right|}$, $\left|\a\right|=\max\{\lambda_{1},\lambda_{3}\}$, and we observe that the visitation times $(a_{1}, a_{2}, a_{3})$ of each sub-interval of $\frac{\a}{\left|\a\right|}$ admit the equation  $a_{2}=a_{1}+a_{3}-1$.
  Consideraing the cyclic approximation of rank 1 stacking structure, we construct a series of cyclic approximation with the base interval to be the second sub-interval of $\a$. Together with the relation $a_{2}=a_{1}+a_{3}-1$, we demonstrate the fundamental structure for whirly property, summarized as Lemma \ref{Major} .\\
   \indent  The last part of the proof is to use a density point argument to extend the fundamental structure based on the Veech Induction to general measurable subsets.
\bcn\label{Conj}
Let $\pi \in \mathcal{G}_{m}^{0}$, $m\geq 3$, for Lebesgue
almost all $\l\in\Lambda_{m}$, the $m$-dimensional cone of positive real numbers, the interval exchange transformation $T_{(\l ,\pi)}$ is whirly.
\ecn

\section{The Space of Three Interval Exchange Transformation}
\indent \indent In W.A.Veech\cite{VEE3}, key results in the theory about interval exchange transformation space
are established. We will utilize the result in W.A.Veech\cite{VEE1} and \cite{VEE3}. Let $m>1$, and specifically here,
let $\pi$ be the symmetric permutation (i.e. $\pi=(m ,m-1,\cdots ,1)$).
In W.A.Veech\cite{VEE2} it is proved that for almost every $\l$ the induced
transformation of $T_{\l ,\pi}$ on $[0, \max
\{\l_{1},\l_{m}\})$ is an $(\a ,\pi)$ interval exchange
transformation with $\left|\a\right|=\max \{\l_{1},
\l_{m}\}$ and $\pi$ still the same symmetric permutation. That is a
transformation $\cT_{2}:(\l,\pi)\rightarrow
(\frac{\a}{\left|\a\right|},\pi)$, or simply, $\cT_{2}(\l) \sim \cT_{2}(\l,\pi)$. So without confusion, let
$\cT_{2}(\l)=\cT_{2}(\l,\pi)$. When $m=3$, $f_{2}(\l)=(\frac{\displaystyle 1}{\displaystyle
1-\l_{1}}+\frac{\displaystyle 1}{\displaystyle
1-\l_{3}})\prod ^{2}_{j=1}\frac{\displaystyle 1}{\displaystyle
\l_{j}+\l_{j+1}}$ is the density of a conservative ergodic
invariant measure for $\cT_{2}$ by W.A.Veech\cite{VEE1}.\\
\indent We claim that if $(\l ,\pi)$ satisfies
i.d.o.c. and $\pi (j)=m-j+1$, there exists some $k$
such that $\cZ^{k}(\l ,\pi)=(\a ,\pi)$ with $\left|\a
\right|=\max \{\l_{1} ,\l_{m}\}$.
To verify this we need the following lemma:
\bl\label{AdmLger}
If $\l\in
\Lambda_{m}, m\geq 3 ,T_{(\l ,\pi)}$ satisfies i.d.o.c. , and
$\cZ^{k}(\l ,\pi)=(\l ' ,\pi ')$, where $k$ is the largest integer such that $\left|\l
'\right|>\max \{\l_{1}, \l_{m}\}$, then $J=[0, \max
\{\l_{1}, \l_{m}\})$ is an admissible interval of
$(\l ' ,\pi ')$.
\el
\bpf
If
$\l_{1}>\l_{m}$, then $\l_{1}$ is a discontinuous
point of $T_{(\l ' ,\pi ')}$, $[0, \l_{1})$ is an
admissible interval of $(\l ' , \pi ')$. If $\l_
{m}>\l_{1}$ , let $\b_{t}'=\b_{t} ' (\l ')= \sum _{i=1}^{t} \l_{i} '$. Since $\l_{m}=T(\b _{m-2})$ and $T_{(\l ' ,\pi ')}$ is the induced transformation of $T_{(\l ,\pi)}$ on $[0, |\l '|)$, we have that there exists $1\leq t\leq m-1$ and a $k_{t}>0$ such that $\l_{m}=T_{(\l ', \pi ')}(\b_{t}')=T_{(\l , \pi )}^{k_{t}}(\b_{t}')$. By the definition of admissible interval,
$[0,\l_{m})$ is an admissible interval associated with
$T_{\l ' ,\pi '}$.
\epf
 \bp
 Suppose $\l \in
\Lambda_{m-1}, \pi (j)=m-j+1$, $1\leq j\leq m$, and $(\l ,\pi)$
satisfies i.d.o.c. Then there exists $k_{0}\in \mN$ such that
$\cZ^{k_{0}}(\l ,\pi)=(\a ,\pi)$, where $\left|\a\right|=\max
\{\l_{1},\l_{m}\}$.
Therefore, $\cT_{2} (\l ,\pi)=\cR^{k_{0}}(\l ,\pi)$.
\ep
\bpf
Assume for all $k\in \mN$, $\cZ^{k}(\l ,\pi)=(\a^{(k)} ,\pi)$
such that $\left|\a^{(k)}\right|\neq \max \{\l_{1},
\l_{m}\}$. Since $(\l ,\pi)$ satisfies i.d.o.c.,
$\left|\pi^{*}_{1}(\cZ^{k}(\l ,\pi))\right|\rightarrow 0$ as
$k\rightarrow \infty$ (see M. Viana\cite{VIA} Corollary 5.2 for a detailed proof), there exist $k_{0}\geq 0$ such that
$\left|\pi^{*}_{1}(\cZ^{k_{0}}(\l ,\pi))\right|>\max
\{\l_{1}, \l_{m}\}$ , and
$\left|\pi^{*}_{1}(\cZ^{k_{0}+1}(\l ,\pi))\right|<\max
\{\l_{1}, \l_{m}\}$. By Lemma \ref{AdmLger} for any
$r>\left|\pi^{*}_{1}(\cZ^{k_{0}+1}(\l, \pi))\right|$, $[0,r)$
is not an admissible interval of $(\l ', \pi ')=\cZ^{k_{0}}(\l ,\pi)$, that
is a contradiction to the fact that $[0, \max
\{\l_{1},\l_{m}\})$ is an admissible interval
 of $(\l ', \pi ')$.
\epf
The above argument assures us that essential general results about the
iteration of Rauzy-Veech induction may be applied to $\cT_{2}$. For
convenience, lets denote the induced map of $T_{\l ,\pi}$ on
$[0,\max \{\l_{1}, \l_{m}\})$ by $(\a ,\pi)$, and define $\cZ_{*}
:\Lambda_{m} \times \{\pi\}\rightarrow \Lambda_{m}\times \{\pi\}$
by $\cZ_{*} (\l ,\pi)=(\a ,\pi)$ with $\left|\a\right|=\max \{\l_{1}, \l_{m}\}$.\\
\indent Next we limit the discussion to the case $m=3$. Recall  Section 1 for the visitation matrix associated with $\cZ^{n}(\l
,\pi)$, $\cZ^{n}(\l ,\pi)=(\a^{(n)}, \pi)$. We have $\l=A^{(n)}\a^{n}$, and the summation of the $i^{th}$ column of $A^{(n)}$, $a^{(n)}_{i}$ is the first return time
 of the $i^{th}$ subinterval of $[0,\left|\a^{(n)}\right|)$
 under $T_{(\l ,\pi)}$. It will be shown that for all
 $n\in \mN$, $a^{(n)}_{2}= a^{(n)}_{1}+ a^{(n)}_{3}-1$.
 In fact we will verify the same equality for a more
 general case. It is done by looking at the Rauzy
 graph for the closed paths based at $\pi =(3,2,1)$.
 The Rauzy class of $\pi =(3,2,1)$ is
 $\{\pi,\pi_{1},\pi_{2}|\pi_{1}=a\pi=(3,1,2),\pi_{2}=b\pi =(2, 3, 1)\}$.\\
$$A(\pi ,a)=A(\pi_{1} ,a)=\left(
\begin{array}{ccc}
1&1&0\\
0&0&1\\
0&1&0
\end{array}\right)$$
$$A(\pi ,b)=A(\pi_{1} ,b)= \left(\begin{array}{ccc}
1&0&0\\
0&1&0\\
1&0&1
\end{array}\right)$$
$$A(\pi_{2} ,a)=\left(\begin{array}{ccc}
1&0&0\\
0&1&1\\
0&0&1
\end{array}\right)$$
$$A(\pi_{2} ,b)=\left(\begin{array}{ccc}
1&0&0\\
0&1&0\\
0&1&1
\end{array}\right) .$$
\bl\label{RecurTime}
If $\cZ_{*}(\l ,\pi)=(\a ,\pi)$,
$\l\in\Lambda_{3}$, $\pi =(3,2,1)$, and the visitation matrix
is $A$ (i.e. $\l =A\a$). Then $a_{2}+1=a_{1}+a_{3}$.
\el
\bpf
To prove this Lemma, we look into the following two cases:
\begin{case}\label{case1}[$ab^{l}a\mbox{ or }ba^{l}b$]
\begin{enumerate}
\item
Starting from $\pi$, go along the path $ab^{l}a$, and come back
to $\pi$. Then the associated visitation matrix is $A^{(l+2)}$, we
want to show that:\\
$$a^{(l+2)}_{2}=a^{(l+2)}_{1}+a^{(l+2)}_{3}-1.$$
Since $$A^{(1)}=A(\pi ,a)=\left( \begin{array}{ccc}
1&1&0\\
0&0&1\\
0&1&0
\end{array}\right)$$
$$a^{(1)}_{1}=a^{(1)}_{3}=1, a^{(1)}_{2}=2$$
$$A^{(2)}=A^{(1)}\cdot A(\pi ,b)=(A^{(1)}_{1},\, A^{(1)}_{2},\, A^{(1)}_{3})\cdot
\left(\begin{array}{ccc}
1&0&0\\
0&1&0\\
1&0&1
\end{array}
\right) $$
$$=(A^{(1)}_{1}+A^{(1)}_{3},\, A^{(1)}_{2},\, A^{(1)}_{3}) ,$$
where $A^{(n)}_{i}$ is the $i$-th column vector of $A^{(n)} .$
$$\cdots\cdots$$
$$A^{(l+1)}=(A^{(1)}_{1}+lA^{(1)}_{3},\, A^{(1)}_{2},\, A^{(1)}_{3})$$
\be
\begin{array}{l}
A^{(l+2)}=A^{(l+1)}\cdot A(\pi_{1} ,a)\\
=(A^{(1)}_{1}+lA^{(1)}_{3} , \, A^{(1)}_{2} , \, A^{(1)}_{3})\cdot\left(
\begin{array}{ccc}
1&1&0\\
0&0&1\\
0&1&0
\end{array}\right)\\
=(A^{(1)}_{1}+lA^{(1)}_{3},\,
A^{(1)}_{1}+(l+1)A^{(1)}_{3},\, A^{(1)}_{2}) \end{array}
\ee
Therefore
\be
\begin{array}{l}
a^{(l+2)}_{1}=a^{(1)}_{1}+la^{(1)}_{3}=l+1\\
a^{(l+2)}_{2}=a^{(1)}_{1}+(l+1)a^{(1)}_{3}=l+2\\
a^{(l+2)}_{3}=a^{(1)}_{2}=2 .
\end{array}
\ee
Thus $a_{2}=a_{1}+a_{3}-1$ is proved for the path $ab^{l}a$.
\item
Similar to the argument in \ref{step01}, if we replace the path of  $ab^{l}a$ to the path of $ba^{l}b$, the associated matrix $A^{(l+2)}$ satisfies:
$$ a^{(l+2)}_{2}=a^{(l+2)}_{1}+a^{(l+2)}_{3}-1 .$$
\end{enumerate}
\end{case}
\begin{case}\label{case2}[$p_{0}ab^{l}a\mbox{ or }p_{0}ba^{l}b$]
\begin{enumerate}
\item\label{step01}
Suppose the closed path is $p=p_{0}ab^{l}a$, where
$p_{0}$ is a closed path based at $\pi =(3,2,1)$, $p_{0}$ admits
length $n_{0}$, and associated with $p_{0}$ is the matrix
$A^{(n_{0})}$ with column summations $a^{(n_{0})}_{1},
a^{(n_{0})}_{2}, a^{(n_{0})}_{3}$ satisfying $ a^{(n_{0})}_{2}+1=
a^{(n_{0})}_{1}+ a^{(n_{0})}_{3}$. Then by similar computation as
Case \ref{case1}. we have the conclusion that, after going along $p$, the
return times satisfy:
\be
a^{(n_{0}+l+2)}_{2}=a^{(n_{0}+l+2)}_{1}+a^{(n_{0}+l+2)}_{3}-1
\ee
\item
Similar to \ref{step01} above, the the same
relation on the three return times is true for the path
$p=p_{0}ba^{l}b.$
\end{enumerate}
\end{case}
\indent By Case \ref{case1} and Case \ref{case2} we have proved Lemma \ref{RecurTime}.
\epf
\section{Whirly Three Interval Exchange Transformations}
\indent \indent Before discussing the 3-interval exchange transformations, let us introduce another way to define the concept of whirly automorphism and verify the equivalence between the two definitions:
\bd[Whirly Automorphism]\label{whirly2}
A rigid ergodic automorphism $T\in G$ is said to be whirly if given $\e >0$, for any $l\in \mN$
(or for any $-l\in \mN$) and a $\mu$-positive measure set $E\in \cB$, there exists $n\in \mN$ such that $T^{n}\in U_{\e}$, and $\mu (T^{n}E\cap T^{l}E)>0$.
\ed
\bt\label{Equiv}
Conditions in
Definition \ref{whirly1} and Definition \ref{whirly2} for an automorphism to be whirly are equivalent to each other.
\et
\bpf
 Suppose $T\in G$ satisfies the condition in Definition \ref{whirly2} (w.l.o.g., we take the case that $-l\in \mN$) ,
 then we claim  that for any $E,F\in \cB$ with $\mu (E),\mu (F)>0$ we have there exists $n\in\mN$ such
 that $T^{n}\in U_{\e}$ and $\mu(T^{n}E\cap F)>0$. Since $T$ is ergodic, there exist $-q\in \mN $
 such that
 $$\mu(T^{q}E\cap F)>0.$$
 \indent Then
 $$\mu(E\cap T^{-q}F)>0 ,$$
 so therefore there exists $n\in \mN$ such that $T^{n}\in U_{\e}$ and
$$\mu (T^{n}(E\cap T^{-q}F)\cap T^{m}(E\cap T^{-q}F))>0.$$
\indent Thus $\mu (T^{n}E\cap F)>0\, .$\\
\indent The opposite direction is obvious.
\epf

\indent Let $\pi$ be the symmetric m-permutation. According to W.A.Veech\cite{VEE3}, there exists $c_{1}, c_{2}, \cdots , c_{n}\in \{a, b\}$, such that: $c_{n}\circ c_{n-1}\circ \cdots \circ c_{1}(\pi )=\pi;$\\
Let $\pi^{(0)}=\pi,\,\pi^{(1)}=c_{1}\pi^{(0)},\,
\pi^{(2)}=c_{2}\pi^{(1)},\cdots ,\pi^{(n)}=c_{n}\pi^{(n-1)}=\pi$, let $A^{(i)}=A(\pi ^{(i-1)},c_{i}),(1\leq i\leq n)$. Then $B=A^{(1)}A^{(2)}\cdots A^{(n)}$ is a positive $m$ $\times$ $m$ matrix.
\br If $\l \in \Lambda_{m}$, then
$\cZ ^{n}(B\l ,\pi)=(\l, \pi)$, and the orbit of
$(B\l ,\pi)$ under $\cZ$
passes the same sequence of permutations $\{\pi^{j},0\leq j\leq n\}$.
\er
Let $\nu(A)=\underset{1\leq i,j,k\leq m}{\max}\{\frac{\displaystyle
a_{ij}}{\displaystyle {a_{ik}}}\}$,
where $A$ is a positive matrix, then:
\be\label{Mv}
a_{i}\leq \nu (A)a_{j},\; 1\leq i,j\leq
m\mbox{ } (a_{i}\mbox{ is the } i\mbox{th column sum of } A)
\ee
\be\label{MvStar}
\nu (MA)\leq \nu (A), \mbox { for any
nonnegative matrix M with at least non zero element.}\;
\ee
We see that $\nu (B)$ and $\nu (B^{t})$ are both positive numbers greater than one. Let $r=\nu (B)$
and $r'=\nu (B^{t})$.\\
\indent Next we fix $m=3$, while still keeping notations as above. Let $l$ be a given positive integer and we will set
up an open set in $\Lambda_{3}\times \{\pi\}$ and do some
computation on the approximation by the Kakutani tower associated with the Rauzy induction.\\
\indent Let $\e_{1} ,\e_{2}$ be two small positive numbers to be specified
for our purpose later. Let $Y^{*}(\e_{1},\e_{2})=\{\a
|\a\in\Lambda_{m},(1-\frac{\e_{1}}{2})\left|\a\right|>\a_{2}>(1-\e_{1})\left|\a\right|
\, and\, $ $(1+\e_{2})\a_{3}>\a_{1}>\a_{3}\}$, an open subset of $\Lambda_{m}$. Let $W(\e_{1},\e_{2})=B^{2}Y^{*}(\e_{1},\e_{2})\times \{\pi\}$, an open subset of $\Lambda_{m}\times \{\pi\}$.\\
Suppose $(\l,\pi)\in \Delta_{2} \times \{\pi\}$, and there
exists $k\in\mN$ such that $\cZ^{k}(\l,\pi)\in
W(\e_{1},\e_{2})$. We know $\xi =B^{2}\a$ for
some $\a \in Y^{*}(\e_{1},\e_{2})$. Then $\l=A^{(k)}\xi$, where $A^{(k)}$ is the visitation matrix
 associated with $\cZ^{k}(\l ,\pi )$,
 $\cZ^{k+2n}(\l ,\pi)=\cZ^{2n}(\xi , \pi)=(\a ,\pi)$, and $\l=A^{(k)}B^{2}\a$.
 Let $A=A^{(k)}B^{2}$. Since $A^{(k)}$ is a non-negative matrix, by \ref{MvStar} we have $\nu(A)\leq \nu (B)=r$.
 Therefore the following arguments may give us a clear view of the stack structure associated with the Veech-induction map $\cT_{2}$:
 \bi
\item[\bf Claim 1\/]  $T^{a_{2}}$ translates the subinterval $I^{\a}_{2}$ (i.e. the
second subinterval of $I^{\a}$) to the left by $(\a_{1}-\a_{3})$. That is
$$ I^{\a}_{2}\cap T^{a_{2}}(I^{\a}_{2})=\left[\a_{1}, \a_{1}+\a_{2}-l(\a_{1}-\a_{3})\right) .$$
Since $l$ is a fixed positive integer, and $\e_{1}$, $\e_{2}$ are small enough, we have
\be\label{Claim1}
\begin{array}{l}
\mu (I^{\a}_{2}\cap T^{la_{2}}(I^{\a}_{2}))=
\a_{2}-l(\a_{1}-\a_{3})>\a_{2}-l\e_{2}\a_{3}\\
>\a_{2}-l\e_{1}\e_{2}\left|\a\right|>\a_{2}-l\frac{\displaystyle {\e_{1}\e_{2}}}{\displaystyle {1-\e_{1}}}\a_{2}\\
=(1-l\frac{\displaystyle{\e_{1}\e_{2}}}{\displaystyle{1-\e_{1}}})\a_{2}\,.
\end{array}
\ee
\item[\bf Claim 2\/]The remainder of the column with base $I^{\a}_{2}$ and
height $a_{2}$ has measure:\\
\be\label{Claim2}
\begin{array}{l} \left|\l \right|-\mu
(\cup^{a_{2}}_{i=0}T^{i}(I^{\a}_{2}))=\\
a_{1}\a_{1} +a_{3} \a_{3}
<r\a_{2}(\a_{1}+\a_{3})
<ra_{2}\e_{1}\left|\a\right|<ra_{2}\frac{\displaystyle{\e_{1}}}{\displaystyle{1-\e_{1}}}\a_{2}
<\frac{\displaystyle{\e_{1}}}{\displaystyle{1-\e_{1}}}\left|\l\right|.
\end{array}
\ee
 From (\ref{Claim1}) and (\ref{Claim2}), for any $\e > 0$, we can select $\e_{1},\, \e_{2}$ small enough such that $T^{la_{2}}\in U_{\e}(Id)$.\\
\item[\bf Claim 3\/]  (the $'$whirly part$'$) $T^{a_{3}}$ sends $[\a_{1}+\a_{2}+(\a_{1}-\a_{3}), \left |
\a\right |)$ to $[\a _{1}- \a_{3},\a _{3})$
which is continuous under $T^{a_{1}}$. That is to say $[\a_{1}+\a_{2}+(\a_{1}-\a_{3}), \left |\a\right |)$ is continuous under $T^{a_{1}+a_{3}}=T^{a_{2}-1}$.
\ei
\indent Similarly by induction:\\
\indent Let
\be
 I^{\a}_{\omega}=[\a_{1}+\a_{2}+l(\a_{1}-\a_{3}),
\left|\a\right|) .
\ee
Then $T^{i}$ are all continuous (linear) on
$I^{\a}_{\o}$ for $i=1,2,\cdots ,l(a_{1}+a_{3})$. And
$T^{l(a_{1}+a_{3})}(I^{\a}_{\omega })\subset I^{\a}_{3}\subset
I^{\a}_{}$.\\
\indent Therefore
$$T^{la_{2}}(I_{\o}^{\a})=T^{l(a_{1}+a_{3}-1)}(I^{\a}_{\o})$$
$$=T^{-l}(T^{l(a_{1}+a_{3})}(I^{\a}_{\o}))\subset T^{-l}(I^{\a}) ,$$
which implies
\be
T^{la_{2}}(I^{\a}_{\o})\subset (T^{la_{2}}(I^{\a}))\cap
T^{-l}(I^{\a}).
\ee
\indent Hence
\be\label{Inter}
\begin{array}{l}\mu (T^{la_{2}}(I^{\a})\cap
T^{-l}(I^{\a}))\\\geq \mu (T^{la_{2}}(I^{\a}_{\o}))=\a_{3}-l(\a_{1}-\a_{3})>\a_{3}-l\e_{2}\a_{3}=(1-l\e_{2})\a_{3} .
\end{array}
\ee
Note: {\bf Claim 1\/} and {\bf Claim 2\/} show that $T^{la_{2}}$ is close to the identity map; \ref{Inter} shows that we are on the right way to the whirly property (Definition \ref{whirly2}).\\
\indent By {\bf Claim 1\/}, {\bf Claim 2\/} and {\bf Claim 3\/}, choosing a positive constant $\mathfrak{C}_{\e
,l}$ associated with $\e$, $l$ and small enough, we have the following Lemma:
\bl\label{Major}
\indent Let $\pi=(3,2,1)$
for almost all $\l\in\Lambda_{3}$, for any
$0<\e <\frac{1}{10}$, $l\in \mN$, there exists $\mathfrak{C}_{\e ,l}$ small
enough such that for $k$ large enough, $\cZ^{k}(\l ,\pi)=(\eta ,\pi)\in
W(\mathfrak{C}_{\e ,l}, \mathfrak{C}_{\e ,l})$.
We have that there exists $n\in \mN$, $\cZ^{k+2n}(\l ,\pi)=(\a ,\pi)$, such that:
$$\begin{array}{l}
P1) \cdots\cdots \mu (I^{\a}\cap
T^{la_{2}}(I^{\a}))>(1-\e)\left|\a\right|\\
P2)\cdots\cdots \left| \l\right|-\mu
(\cup^{a_{2}-1}_{i=0}T^{i}(I^{\a}_{2}))<\e \left|\l\right|\\
P3) \cdots\cdots \mu(T^{la_{2}}(I^{\a})\cap
T^{-l}(I^{\a}))>\frac{\e}{3}\left|\a\right|.\\
\end{array}$$
\el
\indent Now let $N^{(\l)}_{\e ,l}\subset\mN$ be defined by
\be\label{nt}
N^{(\l)}_{\e ,l}=\{n_{t}|n_{1}<n_{2}<\cdots
<n_{i}<\cdots , \cZ^{n_{t}-n}(\l ,\pi)\in W(\mathfrak{C}_{\e ,l},
\mathfrak{C}_{\e ,l})\}.
\ee
\indent By Veech's Ergodic Theorem (W.A.Veech\cite{VEE1} Theorem 1.1) on $\cT_{2}$, we
know that for Lebesgue a.e. $\l \in \Lambda_{3}$, $T(\l ,
\pi)$ is uniquely ergodic (thus ergodic with
 respect to Lebesgue measure), and $N^{\l}_{\e ,l}$ is a set with infinitely many elements,
  for any $0<\e <\frac{1}{10}$, $l\in \mN$. Lets continue to study such $T_{(\l ,\pi)}$. As usual we use $T$ to denote $T_{\l ,\pi}$.\\
\indent We know that
$A=A^{(k)}B^{2}$, $1\leq \nu(A) \leq \nu (B)=r$. We need
 $B^{2}$ here instead of  $B$, in order to get a $T$-stack with the base $I^{\a}$,
 which is a relatively large portion of $I^{\l}$. The following Lemma will be used in the last step, a density point argument, of the proof of Theorem \ref{Main}.
\bl\label{3p5}
All notations as above, let $T=T(\l ,\pi)$, then for $a.e.\, \l \in \Lambda_{3}$ there
exists a positive integer $a_{*}$ such that $T^{i}$ $(1\leq i\leq
a_{*})$ are continuous (linear) on $I^{\a}$, $T^{i}(I^{\a})\cap
T^{j}(I^{\a})=\emptyset$, $(i\neq j, 0\leq i,j <a_{*})$, and
$$a_{*}\left|\a\right|>\frac{1}{b_{M}(1+2r\cdot r')}\left|\l\right| ,$$
where $b_{M}=\max \{b_{11},b_{12},b_{13}\}$.
\el
\begin{proof}
 \indent We know that $A=A^{(k)}B^{2}$.\\
\indent Suppose $\cZ^{k+n}(\l ,\pi)=(\eta ,\pi)=(B\a ,\pi)$, where $\eta=(\eta_{1},\eta_{2},\eta_{3})$, $I^{\eta}_{1}=[0,
\eta_{1})$, $\eta_{1}=b_{11}\a_{1}+b_{12}\a_{2}+b_{13}\a_{3}$, and
\be\label{Eta1} I^{\a}\subset I^{\eta}_{1} .\ee
\indent Meanwhile $\eta_{1}<b_{M}(\a_{1}+\a_{2}+\a_{3})<b_{M}\left|\a\right|$; that
is
\be\label{AlEta}
\left|\a\right|> \frac{1}{b_{M}}\eta_{1} .
\ee
\indent At the same time, since
$$\eta_{1}=b_{11}\a_{1}+b_{12}\a_{2}+b_{13}\a_{3}$$
$$\eta_{2}=b_{21}\a_{1}+b_{22}\a_{2}+b_{23}\a_{3}$$
$$\eta_{3}=b_{31}\a_{1}+b_{32}\a_{2}+b_{33}\a_{3} ,$$
it follows that
\be\label{Eta}
\eta_{2}, \eta_{3}<r'\eta_{1} .
\ee
\indent Remembering that $\l =A^{(k)}B\eta$,  i.e.
$\l=a^{(k+n)}_{1}\eta_{1}+a^{(k+n)}_{2}\eta_{2}+a^{(k+n)}_{3}\eta_{3}$,
by \ref{Mv} and \ref{MvStar} we have
$$a^{(k+n)}_{2}, a^{(k+n)}_{3}<r a^{(k+n)}_{1} ,$$
and by \ref{Eta} we have
\be\label{Cover}
a^{(k+n)}_{1}\eta_{1}>\frac{1}{1+2rr'}\left|\l\right| .
\ee
\ref{AlEta} and \ref{Cover} imply that $
a^{(k+n)}_{1}\left|\a\right|>\frac{1}{b_{M}(1+2rr')}\left|\l\right| $, combining this with
\ref{Eta1}, the Lemma is proved.
\end{proof}
{\em\bf Proof of Theorem \ref{Main} (with a Density Point Argument)\/\/}
\bpf
\indent Let $\pi=(3,2,1)$, and let $\l$ be in the full measure subset of $\Lambda_{3}$ as required by Lemma \ref{Major}.\\
\indent Define
$\mathfrak{G}=\underset{N=1}{\overset{\infty}{\cap}}\underset{\overset{t\leq
N}{n_{t}\in N^{(\l)}_{\e ,l}}}{\cup}\mathfrak{G}_{t}$,
where $\mathfrak{G}_{t}=\cup ^{a_{*}^{(n_{t})}-1}_{i=0}T^{i}(I^{\a^{n_{t}}})$, with $n_{t}$ as defined in \ref{nt}.\\
\indent According to Lemma \ref{3p5}, $\mu (\mathfrak{G})\geq \frac{\displaystyle 1}{\displaystyle b_{M}}\frac{\displaystyle 1}{\displaystyle 1+2rr'}\left|\l\right|$. \\
\indent Suppose $E$ is an arbitrary measurable set, $E\subset
[0,\left|\l\right|)$, $\mu (E)>0$. Then by the ergodicity of $T$ there exists $q\in \mN$
such that $\mu (T^{-q}(E)\cap \mathfrak{G})>0$. Therefore, by the Lebesgue Density Theorem, there exists a
point of density one $x\in T^{-q}(E)\cap \mathfrak{G}$. By definition of $\mathfrak{G}$, we have $x\in T^{-q}_{(\l ,\pi)}(E)\cap J_{k}$, where the left close right open interval
$J_{k}=T^{i_{k}}(I^{\a^{(S_{k})}})$, $S_{k}=n_{t_{k}}$
, $0\leq i_{k}<a^{(S_{k})}_{*}$, and the approximate density satisfies:
\be\label{Density}
\underset{k\rightarrow \infty}{\lim}\frac{\mu (T^{-q} (E)\cap
J_{k})}{\mu (J_{k})}=1 .
\ee
\indent By Lemma \ref{Major} , we know that since
$S_{k}=n_{t_{k}}\in
N^{(\l)}_{\e ,l}$.
\be\label{Lap}
\begin{array}{l}
\mu ((T^{la^{(S_{k})}_{2}}J_{k})\cap T^{-l}(J_{k}))\\
\cdots\cdots =\mu
(T^{i_{k}}(T^{la^{(S_{k})}_{2}}(I^{\a^{(S_{k})}})\cap T^{-l}(I^{\a
^{(S_{k})}})))\\
$$>\frac{\e}{3}\left|\a^{(S_{k})}\right|.
\end{array}
\ee
\indent \ref{Density} implies there exists $k_{0}$ such that
$$\frac{\mu (T^{-q}(E)\cap J_{k_{0}})}{\mu (J_{k_{0}})}>(1-\frac{\e}{10}) .$$
\indent Therefore by \ref{Lap} we have
$$\mu (T^{la^{(S_{k})}_{2}}(T^{-q} (E))\cap (T^{-q} (E))>0 .$$
\indent Thus $\mu (T^{la^{(S_{k})}_{2}}(E)\cap T^{-l} (E))>0$. Since $S_{k}=n_{t_{k}}\in N^{(\l)}_{\e ,l}$, together with Theorem \ref{Equiv} and Lemma \ref{Major}, we have
proved Theorem \ref{Main}.
\epf
\bc
Let $\pi=(3,2,1)$, for Lebesgue almost all $\l\in
\Lambda_{3}$, the interval exchange transformation $(\mX ,\cB ,T_{(\l, \pi)})$ admits no nontrivial spatial factor.
\ec
\bpf
By  Proposition 1.9 of E.Glasner, B.Weiss\cite{GW}.
\epf
\renewcommand{\abstractname}{Comments and Acknowledgements}
\begin{abstract}
 The result about 3-interval exchange transformations is a part of the author Y. Wu's 2006 Ph.D. Thesis at Rice University, Department of Mathematics.
 The author Y. Wu thanks Rice University for posting his Doctor of Philosophy Thesis online at the Rice University Digital Scholarship Archive. He would also like to thank his Ph.D. advisor W. A. Veech for directing his Ph.D. Thesis.
\end{abstract}

\end{document}